\documentclass[final,1p,times]{elsarticle}

\usepackage{amsmath,amssymb}
\usepackage[T1]{fontenc}
\usepackage[cp1250]{inputenc}

\newcommand{\nat}{\mbox{$\mathbb{N}$}}
\newcommand{\real}{\mbox{$\mathbb{R}$}}

\newcommand{\Qb}{\mbox{$\mathbb{Q}$}}
\newcommand{\Mcl}{\mbox{$\mathcal{M}$}}
\newcommand{\fcont}{f_{\textit{cont}}}
\newcommand{\fsing}{f_{\textit{sing}}}
\newcommand{\fpp}{f_{\textit{pp}}}
\newcommand{\gcont}{g_{\textit{cont}}}
\newcommand{\gsing}{g_{\textit{sing}}}
\newcommand{\gpp}{g_{\textit{pp}}}
\newcommand{\sn}{\succeq_n}

\newtheorem{theorem}{Theorem}[section]
 \newtheorem{corollary}[theorem]{Corollary}
 \newtheorem{lemma}[theorem]{Lemma}
 \newtheorem{proposition}[theorem]{Proposition}
 \newtheorem{observation}[theorem]{Observation}
 
 \newtheorem{notation}[theorem]{Notation}
 \newdefinition{remark}[theorem]{Remark}
 \newdefinition{example}[theorem]{Example}

 \newproof{proof}{Proof}
 \newproof{pfs}{Proof of Theorem~\ref{th:support}}
 \newproof{pfd}{Proof of Lemma~\ref{lm:delta}}

\numberwithin{equation}{section}

\journal{Journal of Mathematical Analysis and Applications}

\begin{document}

\begin{frontmatter}



\title{New integral representations of $n$th order convex functions}


\author{Teresa Rajba}
\ead{trajba@ath.bielsko.pl}
\address{Department of Mathematics and Computer Science, University of Bielsko-Biała, ul. Willowa 2, 43-309 Bielsko-Biała, Poland}

\begin{abstract}
In this paper we give an integral representation of an $n$-convex function $f$ in general case without additional assumptions on function $f$. We prove that any $n$-convex function can be represented as a sum of two $(n+1)$-times monotone functions and a polynomial of degree at most $n$. We obtain a decomposition of $n$-Wright-convex functions which generalizes and complements results of Maksa and Pales \cite{MaksaPales2009}. We define and study relative $n$-convexity of $n$-convex functions. We introduce a measure of $n$-convexity of $f$. We give a characterization of relative $n$-convexity in terms of this measure, as well as in terms of $n$th order distributional derivatives and Radon-Nikodym derivatives. We define, study and give a characterization of strong $n$-convexity of an $n$-convex function $f$ in terms of its derivative $f^{(n+1)}(x)$ (which exists  a.e.) without additional assumptions on differentiability
of $f$.
We prove that for any two $n$-convex functions $f$ and $g$, such that $f$ is $n$-convex with respect to $g$, the function $g$ is the support for the function $f$ in the sense introduced by Wasowicz \cite{Wasowicz2007}, up to polynomial of degree at most $n$.
\end{abstract}


\end{frontmatter}

% \linenumbers

\section{Introduction}

The notion of \textit{$n$th order convexity} (or \textit{$n$-convexity}) was defined in terms of divided differences by Popoviciu \cite{Popoviciu1934} (cf. also \cite{RobertsVarberg1973}, \cite{Kuczma1985}), however, we will not state it here. Instead we list some definitions of $n$th order convexity which are equivalent to the Popoviciu's definition.
\begin{proposition}\label{prop:1.1}
A function $f(x)$ is $n$-convex on $(a,b)$ $(n \geqslant 1)$ if and only if its derivative $f^{(n-1)}(x)$  exists and is convex on $(a,b)$ (with the convention $f^{(0)}(x)=f(x)$).
\end{proposition}
This fact first was proved by Hopf \cite[p. 24]{Hopf1926} and by Popoviciu \cite[p. 38]{Popoviciu1934} (see also \cite{Kuczma1985}, \cite{RobertsVarberg1973}). Many results on $n$-convex functions one can found, among others, in \cite{KarlinStudden1966}, \cite{BessenyeiPales2003}, \cite{Granata1982}, \cite{Kuczma1985}, \cite{RobertsVarberg1973}, \cite{NikodemPales2007}, \cite{PinkusWulbert2005}, \cite{Wasowicz2006}, \cite{Wasowicz2007}, \cite{Wasowicz2010}, \cite{GilanyiPales2008}.

Recall that convex functions satisfy various smoothness properties. A convex function defined on $(a,b)$ is continuous and has both right and left derivatives $f'_R(x)$ and $f'_L(x)$ at each point of $(a,b)$. In addition both these derivatives are non-decreasing and satisfy inequality $f_L'(x) \leqslant f_R'(x)$ for all $x \in (a,b)$ (see \cite{RobertsVarberg1973}, \cite{Kuczma1985}). Thus we have
\begin{proposition}\label{prop:1.2}
A function $f(x) \colon (a,b) \to \real$ is $n$th order convex $(n\geqslant 1)$ if and only if its right derivative $f_R^{(n)}(x)$ (or left derivative $f_L^{(n)}(x)$) exists and is non-decreasing on $(a,b)$.
\end{proposition}
If $f(x)$ is sufficiently smooth on $[a,b]$, then from Taylor's Theorem we have
$$
f(x) = \sum_{k=0}^{n} \frac{ f^{(k)}(a) (x-a)^k }{ k! } + \frac{1}{n!} \int_{a}^b (x-t)_+^{n} f^{(n)}(t) dt,
$$
where $(x-t)^{n-1}_+ = \max \{ (x-t)^{n-1}, 0 \}$.

Now assume $f(x)$ is $n$th order convex on $(a,b)$ $(n \geqslant 1 )$. Then the left and right derivatives $f_L^{(n)}(x)$ and $f_R^{(n)}(x)$ exist on $(a,b)$. In addition, both these functions are non-decreasing. With each such $f$ we associate the measure $\mu$ defined on $(a,b)$ by
$$
\mu( [x,y] ) = f_R^{(n)}(y) - f_L^{(n)}(x),
$$
for $a < x \leqslant y < b$. This is a non-negative Borel measure on $(a,b)$. If $f_R^{(n)}(a)$ is finite then $\mu$ can be extended to a bounded (finite) measure on whole $[a,c]$, for all $c<b$. In this case $f(x)$ has the representation
$$
f(x) = \sum_{k=0}^{n} \frac{ f_R^{(k)}(a) (x-a)^k }{ k! } + \frac{ 1 }{ n! } \int_a^b (x-t)^{n}_+ d \mu (t),
$$
for $x \in (a,b)$. If we cannot extend $\mu$ to the endpoint $a$, then we will have this representation only on closed subintervals of $(a,b)$. The converse also holds. These results can be found in Popoviciu \cite{Popoviciu1934} (see also Karlin and Studden \cite{KarlinStudden1966}, Bullen \cite{Bullen1971}, Brown \cite{Brown1989}, Granata \cite{Granata1982}, Pinkus and Wulbert \cite{PinkusWulbert2005}). In other words, the above integral representation is valid for all $x \in (a,b)$ if $\mu$ is of bounded variation on $(a,b)$, otherwise we have this representation only on closed subinterval of $(a,b)$.

In this paper we give an analogue of the integral representation above in general case. The representation we obtain deals with measures $\mu$ with not necessarily bounded variations. Our characterization is constructive. We give explicit formulas for $n$-spectral measures corresponding to an $n$-convex function in this representation (see Section 2).

The strength of the representation developed in Section 2 is exploited in the rest of the paper. It is used to further study of $n$-convexity, and to obtain complete characterizations of strong $n$-convexity, $n$-Wright-convexity, and  relative $n$-convexity of functions, among other. Finally, the representation is employed to examine support-type properties of $n$-convex functions. 

In Section 3 we prove that an $n$-convex function can be represented as a sum of two $(n+1)$-times monotone functions and a polynomial of degree at most $n$. This result generalizes the well-known theorem on representation of a convex function as a sum of non-increasing and non-decreasing functions, and a polynomial of degree at most 1 (see Roberts and Varberg \cite{RobertsVarberg1973}). Using our decomposition we obtain the decomposition of $n$-Wright-convex functions, which generalizes and complements results of Maksa and P\'{a}les \cite{MaksaPales2009}.

In Section 4 we define and study relative $n$-convexity of $n$-convex functions. Relative $n$-convexity induces the partial ordering in the set of $n$-convex functions. We define a measure of $n$-convexity of an $n$-convex function $f$ using $n$-spectral measures in our representation. We give a characterization of relative $n$-convexity in terms of the measure of $n$-convexity, as well as in terms of $n$th order distributional derivatives, and in terms of Radon-Nikodym derivatives. Using the Lebesgue decomposition of $n$-spectral  measures corresponding to an $n$-convex function $f$, we consider the corresponding decomposition of the function $f$. This decomposition is applied to derive some useful characterizations of the relative $n$-convexity.

We define and study the notion of strong $n$-convexity that generalizes the strong convexity. It is well known that the strong convexity of a function $f$ can be characterized in terms of its second derivative $f''(x)$ for twice differentiable $f$.
We give a characterization of the strong $n$-convexity of an $n$-convex function $f$ in terms of only derivative $f^{(n+1)}(x)$ (which exists almost everywhere with respect to Lebesgue measure), without any additional assumptions on differentiability of $f$.

In Section 5 we obtain a generalization of Wasowicz \cite{Wasowicz2007} results. We prove, that for any two $n$-convex functions $f$ and $g$, such that $f$ is $n$-convex with respect to $g$, the function $g$ is the support for the function $f$ in the sense introduced by Wasowicz \cite{Wasowicz2007}, up to a polynomial of degree at most $n$.

\section{Integral representation}

In this chapter we give an integral representation of an $n$-convex function $f$ without additional assumptions on $f$. We derive explicit formulas for $n$-spectral measures corresponding to $f$ that can be applied to measures of not necessary bounded variation on $(a,b)$.

By $\lambda$ we denote the Lebesgue measure.
Let $\Pi_n$ be the family of all polynomials of degree at most $n$. Let $f \colon (a,b) \to \real$ be an $n$th order convex function on the interval $(a,b)$, where $- \infty \leqslant a < b \leqslant \infty$, $n=1,2,\ldots$. Then $f_R^{(n)}(x)$ is non-decreasing and right-continuous on $(a,b)$. Henceforth $f^{(n)}(x)$ will be used to denote $f^{(n)}_R(x)$. A function $f^{(n)}(x)$ must satisfy one of the following three conditions:
\begin{itemize}
	\item[A.] There exist $x_1,x_2 \in (a,b)$ such that $f^{(n)}(x_1) < 0$ and $f^{(n)}( x_2 ) >0$,
	\item[B.] $f^{(n)}(x) \geqslant 0$ for all $x \in (a,b)$,
	\item[C.] $f^{(n)}(x) \leqslant 0$ for all $x \in (a,b)$.
\end{itemize}

\begin{theorem}\label{th:1.1}
For $n \geqslant 1$ each $n$th order convex function $f \colon (a,b) \to \real$ satisfying the property A admits the representation of the form
\begin{equation}\label{eq:1.1}
f(x) = \int_{(a,\xi]} (-1)^{n+1} \frac{ [-(x-u)]^n_+ }{ n! } dg_{(n)-}(u) + \int_{[\xi,b)} \frac{(x-u)^n_+}{n!} dg_{(n)+}(u) + Q(x),
\end{equation}
where $\xi \in (a,b)$, $g_{(n)-} \leqslant 0$ is a non-decreasing right-continuous function on $(a,b)$, $g_{(n)+}$ is a non-decreasing left-continuous function on $(a,b)$ such that $g_{(n)+}(\xi) = g_{(n)-}(\xi) = 0$, and $Q \in \Pi_{n-1}$. Moreover, the functions $g_{(n)-}$, $g_{(n)+}$ and $Q$ are determined uniquely, $g_{(n)+} = f_+^{(n)}$ a.e., $g_{(n)-} = f_-^{(n)}$ a.e.
% $Q(x) = f(x) - \Psi_n( x, (a,b), \xi, dg_{(n)-}(u),dg_{(n)+}(u) )$ ($\Psi_n$ is defined by the formula (\ref{eq:1.4})).
\end{theorem}

\begin{notation}\label{not:2}
The quantities
\begin{equation}\label{eq:1.2}
\Psi_{(n)-}(x) =
\Psi_{(n)-}( x; a, \xi, dg_{(n)-}(u) ) = 
\int_{(a,\xi]} (-1)^{n+1} \frac{ [-(x-u)]^n_+ }{ n! } dg_{(n)-}(u),
\end{equation}

\begin{equation}\label{eq:1.3}
\Psi_{(n)+}(x) =
\Psi_{(n)+}( x; \xi, b, d g_{(n)+}(u) )= 
\int_{[\xi,b)} \frac{(x-u)^n_+}{n!} d g_{(n)+}(u)
\end{equation}
appear frequently and hence from now we will be using the above notation.
\end{notation}

\begin{remark}\label{remark:3}
A straightforward calculation shows that
\begin{eqnarray}\label{eq:1.5}
\frac{ d^n }{ dx^n } \Psi_{(n)-}( x; a, \xi, d g_{(n)-}(u) ) & = & \int_{(a,\xi]} [-\chi_{(-\infty,0)}(x-u)] d g_{(n)-}(u) \nonumber \\
& = & g_{(n)-}(x) \textit{ a.e. } ( x \in (a,\xi) ),
\end{eqnarray}
\begin{eqnarray}\label{eq:1.6}
\frac{ d^n }{ dx^n } \Psi_{(n)+}( x; \xi, b, g_{(n)+}(u) ) & = & \int_{[\xi,b)]} \chi_{(0,\infty)}(x-u)d g_{(n)+}(u) \nonumber \\
& = & g_{(n)+}(x) \textit{ a.e. } ( x \in (\xi,b) ).
\end{eqnarray}
\end{remark}

\begin{proof}[Proof of Theorem \ref{th:1.1}]
Let $f$ be an $n$th order convex function satisfying the property A. Then there exists $\xi \in (a,b)$ such that $f^{(n)}(\xi+) \geqslant 0$ and $f^{(n)}(\xi-) \leqslant 0$. Let $g_{(n)-}(x)$ and $g_{(n)+}(x)$ be right-continuous and left-continuous functions, respectively, and such that
\begin{equation}\label{eq:1.8}
g_{(n)-}(x) = \min \{ 0, f^{(n)}(x) \}, g_{(n)+}(x) = \max \{ 0, f^{(n)}(x) \} \textit{ a.e. }
\end{equation}
Then $g_{(n)-}(\xi) = g_{(n)+}(\xi) = 0$. From (\ref{eq:1.5}), (\ref{eq:1.6}) and (\ref{eq:1.8}) we obtain that the functions $f(x)$ and $\Psi_{(n)-}(x)+\Psi_{(n)+}(x)$ differ on $(a,b)$ by a polynomial of degree at most $n-1$. Thus (\ref{th:1.1}) is satisfied. Conversely, assume $f$ is of the form (\ref{eq:1.1}). By Remark \ref{remark:3}, $f_-^{(n)}(x) = g_{(n)-}(x)$ and $f_+^{(n)}(x) = g_{(n)+}(x)$ a.e. Thus $f^{(n)}(x)$ is non-decreasing and right-continuous on $(a,b)$. This implies that $f(x)$ is $n$th order convex on $(a,b)$. The proof is completed.
\end{proof}

\begin{remark}\label{remark:2.4}
Note that since $[-(x-u)]^n_+ = 0$ for $u < x$, and $(x-u)^n_+ = 0$ for $u>x$, the integral (\ref{eq:1.2}) is over $[x,\xi]$ and the integral (\ref{eq:1.3}) is over $[\xi,x]$. Since $d g_{(n)-}(u)$ and $d g_{(n)+}(u)$ are of bounded variations on $(x,\xi)$ and $(\xi,x)$, respectively, the integrals in (\ref{eq:1.2}) and (\ref{eq:1.3}) are well defined.
\end{remark}

\begin{remark}\label{remark:2.5}
If $g_{(n)-}(b-) = 0$, then in (\ref{eq:1.5}) we set $\xi = b$. Similarly if $g_{(n)+}(a+) = 0$, then we put $\xi = a$ in (\ref{eq:1.6}).
%In the case $\xi = b$, the formula (\ref{eq:1.5}) applies for a non-positive non-decreasing function $g_{(n)-}(x)$ satisfying $g_{(n)-}(b-) = 0$. In the case $\xi = a$ the formula (\ref{eq:1.6}) applies for a non-negative non-decreasing function $g_{(n)+}(x)$ satisfying $g_{(n)+}(a+) = 0$.
\end{remark}

\begin{theorem}\label{th:2.6}
For $n \geqslant 1$ each $n$-convex function $f \colon (a,b) \to \real$ satisfying the property $B$ admits the representation
\begin{equation}
f(x) = \int_a^b \frac{(x-u)_+^n}{n!}d g_{(n)}(u) + Q(x),
\label{eq:2.9}
\end{equation}
where $Q(x) = c_n x^n/n!+\ldots+c_0$, $c_n \geqslant 0$, and $g_{(n)}(x)$ is a non-negative non-decreasing left-continuous function on $(a,b)$ satisfying $g_{(n)}(a+) = 0$. Moreover, $Q(x)$ and $g_{(n)}(x)$ are uniquely determined,
$c_n = f^{(n)}(a+)$, 
$g_{(n)}(x) = f^{(n)}(x) - c_n$ a.e., and
$Q(x) = f(x) - \psi_{(n)+}( x; a,b, dg_{(n)}(u))$.
\end{theorem}

\begin{proof}
Assume $f$ is $n$-convex function such that $f^{(n)}(x) \geqslant 0$ $(x \in (a,b) )$. Taking into account that $f^{(n)}(x)$ is non-negative and non-decreasing on $(a,b)$, $f^{(n)}(a+) = c_n$ exists and is finite. Let $g_{(n)}(x)$ be a left-continuous function such that $g_{(n)}(x) = f^{(n)}(x) - c_n$ a.e. $(x \in (a,b))$. Then $g_{(n)}(x)$ is non-negative, non-decreasing, and satisfies $g_{(n)}(a+) = 0$. In view of Remark \ref{remark:2.5}, by (\ref{eq:1.6}) with $a$ in place of $\xi$ and $g_{(n)}(x)$ in place of $g_{(n)+}(x)$, we have
$$
\frac{d^n}{dx^n} \psi_{(n)+}(x; a,b,dg_{(n)}(u)) = g_{(n)}(x) \quad a.e. \quad (x \in (a,b) ).
$$
Consequently
$$
\frac{d^n}{dx^n} \psi_{(n)+}(x; a,b,dg_{(n)}(u)) = f^{(n)}(x)-c_n \quad a.e. \quad (x \in (a,b)).
$$
Thus the functions 
$\psi_{(n)+}(x; a,b, dg_{(n)}(u))$ and 
$f(x) - c_n x^n/n!$ differ on $(a,b)$ by a polynomial of degree at most $(n-1)$. The theorem is proved.
\end{proof}

\begin{theorem}\label{th:2.7}
For $n \geqslant 1$ each $n$th order convex function $f \colon (a,b) \to \real$ satisfying the property $C$ admits the representation of the form
\begin{equation}\label{eq:2.10}
f(x) = \int_a^b (-1)^{n+1} \frac{[-(x-u)]^n_+}{n!} d g_{(n)}(u) + Q(x),
\end{equation}
where $Q(x) = c_n x^n/n! + \ldots + c_0$, $c_n \leqslant 0$, and $g_{(n)}(x)$ is a non-positive non-decreasing right-continuous function on $(a,b)$ such that $g_{(n)}(b-)=c_n$. Moreover, $g_{(n)}$ and $Q$ are uniquely determined, $c_n=f^{(n)}(b-)$, $g_{(n)}(x) = f^{(n)}(x) - c_n$ a.e., and $Q(x) = f(x) - \psi_{(n)-}(x;a,b,d g_{(n)}(u))$.
\end{theorem}

\begin{proof}
The proof is similar to the proof of Theorem \ref{th:2.6} and hence it is omitted.
\end{proof}

\begin{remark}\label{remark:2.8}
The representations (\ref{eq:1.1}), (\ref{eq:2.9}) and (\ref{eq:2.10}) can be rewritten in equivalent forms using the following two measures associated with the distribution functions $g_{(n)-}(x)$ and $g_{(n)+}(x)$, defined as
$$
\begin{array}{lll}
\mu_{(n)-}(du) & = & dg_{(n)-}(u), \\
\mu_{(n)+}(du) & = & dg_{(n)+}(u).
\end{array}
$$
We will call $\mu_{(n)-}$ and $\mu_{(n)+}$ the \textit{$n$-spectral measures} 
of an $n$-convex function $f$.
\end{remark}

The following theorem summarizes Theorems \ref{th:1.1}, \ref{th:2.6} and \ref{th:2.7}.

\begin{theorem}\label{th:2.9}

\noindent
\begin{enumerate}[\upshape a)] 
	\item For $n \geqslant 1$ each $n$-convex function $f \colon (a,b) \to \real$ admits the representation
\begin{equation}
f(x) = \int_{(a,\xi]} (-1)^{n+1} \frac{ [-(x-u)]^n_+ }{ n! } \mu_{(n)-}(du)
 + \int_{[\xi,b)} \frac{(x-u)_+^n}{n!} \mu_{(n)+}(du) + Q(x),
\label{eq:2.11}
\end{equation}
where $\xi \in [a,b]$. Moreover, if $f^{(n)}(x)$ satisfies the condition $B$ (or $C$), then $\xi = a$, $\mu_{(n)-} = 0$ and $\mu_{(n)+}(du) = d(f^{(n)}(u) - f^{(n)}(a+))$
(or $\xi = b$, $\mu_{(n)+} = 0$ and $\mu_{(n)-}(du) = d(f^{(n)}(u) - f^{(n)}(b-))$), and if $f^{(n)}(x)$ satisfies the condition $A$, then $\xi \in (a,b)$, $f^{(n)}(\xi-) \leqslant 0$, $f^{(n)}(\xi+) \geqslant 0$, $\mu_{(n)-}(du) = d f_{-}^{(n)}(u)$, $\mu_{(n)+}(du) = d f_+^{(n)}(u)$ and $Q(x) \in \Pi_n$.

\item If $f^{(n)}(a+) = \alpha$ exists and is finite, then $f(x)$ can be rewritten in the form
$$ %\begin{equation}
f(x) = \int_a^b \frac{(x-u)^n_+}{n!} \mu_{(n)a+}(du) + Q_a(x),
%\label{eq:2.13}
$$ %\end{equation}
where $\mu_{(n)a+}(du) = d(f^{(n)}(u)-\alpha)_+$, $Q_a(x) \in \Pi_n$.

\item If $f^{(n)}(b-) = \beta$ exists and is finite, then $f(x)$ can be rewritten in the form
$$ %\begin{equation}
f(x) = \int_a^b (-1)^{n+1} \frac{[-(x-u)]_+^n}{n!} \mu_{(n)b-}(du) + Q_b(x),
%\label{eq:2.14}
$$ %\end{equation}
where $\mu_{(n)b-}(du) = d(f^{(n)}(u)-\beta)_-$, $Q_b(x) \in \Pi_n$.
\end{enumerate}
\end{theorem}

Denoting
\begin{eqnarray*}
&& \psi_f (x) = \psi_{(n)f}(x) = \\
&& \ \quad \int_{(a,\xi]} (-1)^{n+1} \frac{ [-(x-u)]^n_+ }{ n! } \mu_{(n)-}(du)
 + \int_{[\xi,b)} \frac{(x-u)_+^n}{n!} \mu_{(n)+}(du),
%\label{eq:2.12}
\end{eqnarray*}
(\ref{eq:2.11}) can be rewritten in the form
$$ %\begin{equation}
f(x) = \psi_f(x) + Q(x).
%\label{eq:2.13a}
$$ %\end{equation}

Note, that every function $f(x)$ can be trivially written as $f(x) = f(x) - c x^n/n! + c x^n/n!$ $(c \in \real)$. Thus $f(x)$ can be also written in the form
\begin{equation}
f(x) = \psi_{f-c x^n/n!}(x) + Q_c(x),
\label{eq:2.14a}
\end{equation}
where $Q_c(x) \in \Pi_n$.

Another representation is given in the following theorem. This representation is important in applications of the theory to study relative $n$-convexity.

\begin{theorem}\label{th:2.10}
Let $f \colon (a,b) \to \real$ be an $n$-convex function. For every $\xi \in (a,b)$ the function $f(x)$ has the representation
\begin{equation}
f(x) = \int_{(a,\xi]} (-1)^{n+1} \frac{[-(x-u)]^n_+}{n!} \mu_{(n)\xi-}(du)+
 \int_{[\xi,b)} \frac{(x-u)^n_+}{n!} \mu_{(n)\xi+}(du)+Q_\xi(x),
 \label{eq:2.15}
\end{equation}
where
\begin{eqnarray}
\mu_{(n)\xi-}(du) & = & d[f^{(n)}(u)-f^{(n)}(\xi+)]_-, \nonumber \\
\mu_{(n)\xi+}(du) & = & d[f^{(n)}(u)-f^{(n)}(\xi+)]_+, \nonumber \\ 
Q_\xi \in \Pi_n.  &&
\label{eq:2.16}
\end{eqnarray}
Moreover, we have
\begin{equation}
\mu_{(n)\xi-}+\mu_{(n)\xi+} = \mu_{(n)-}+\mu_{(n)+},
\label{eq:2.17}
\end{equation}
where $\mu_{(n)-}$ and $\mu_{(n)+}$ are the $n$-spectral measures corresponding to $f$.
\end{theorem}

\begin{proof}
Let $a < \xi < b$. Put $c = f^{(n)}(\xi+)$ and denote $g_c(x) = f(x) - c x^n/n!$. Then $g_c^{(n)}(\xi-) \leqslant 0$ and $g_c^{(n)}(\xi+) \geqslant 0$, and consequently the function $g_c^{(n)}(x)$ satisfies the condition $A$. By Theorem \ref{th:2.9} with $g_c(x)$ in place of $f(x)$, and taking into account (\ref{eq:2.14a}), we obtain the representation (\ref{eq:2.15}) with the measures $\mu_{(n)\xi-}$ and $\mu_{(n)\xi+}$ satisfying (\ref{eq:2.16}). It implies that $f^{(n)}(x)-f^{(n)}(\xi+)$ is the distribution function corresponding to the sum of measures $\mu_{(n)\xi-}+\mu_{(n)\xi+}$. By Theorem \ref{th:2.9}, the distribution function corresponding to the sum of $n$-spectral measures $\mu_{(n)-}+\mu_{(n)+}$ equals $f^{(n)}(x)$ up to a constant. Thus the distribution functions of measures $\mu_{(n)\xi-}+\mu_{(n)\xi+}$ and $\mu_{(n)-}+\mu_{(n)+}$ differ on $(a,b)$ by a constant. Consequently these measures coincide, so (\ref{eq:2.17}) is proved.
\end{proof}

\section{$n$-convexity and multiple monotonicity}

From Theorem \ref{th:2.9} on the representation of an $n$-convex function $f$ we obtain that $f$ can be represented by the sum of two $(n+1)$-times monotone functions and a polynomial of degree at most $n$. Applying this we obtain a theorem on decomposition of an $n$-Wright-convex function, which complements and generalizes results of Maksa and P\'{a}les \cite{MaksaPales2009}.

By the standard definition (cf. Williamson \cite{Williamson1956}) a function $f \colon (a,b) \to \real$ is called \textit{$n$-times monotone non-increasing} $(n \geqslant 2)$ if $(-1)^k f^{(k)}(x)$ is non-negative, non-increasing, and convex for $x \in (a,b)$ and $k=0,1,\ldots,n-2$. When $n=1$, $f(x)$ is simply non-negative and non-increasing. The well-known representation for $n$-times monotone non-increasing functions on $(0,\infty)$ states that
\begin{equation}\label{eq:3.1}
f(x) = \int_0^\infty (1-ux)^{n-1}_+ d \beta(u) \quad (x>0),
\end{equation}
with $\beta(u)$ being non-decreasing (see Williamson \cite{Williamson1956}).

A function $f$ is called \textit{$n$-times monotone non-decreasing} (briefly \textit{$n$-times monotone}) $(n\geqslant 2)$ if $f^{(k)}(x)$ is non-negative, non-decreasing, and convex for $x \in (a,b)$ and $k=0,1,2,\ldots,n-2$. When $n=1$, $f(x)$ is simply non-negative and non-decreasing. From (\ref{eq:3.1}) we derive the following representations of functions $f \colon (a,b) \to \real$, which are $(n+1)$-times monotone non-increasing and $(n+1)$-times monotone non-decreasing on $(a,b)$, respectively:
\begin{equation}\label{eq:3.4}
f(x) = \int_a^b \frac{(x-u)^n_+}{n!} d \beta(u),
\end{equation}
\begin{equation}\label{eq:3.5}
f(x) = \int_a^b \frac{ [-(x-u)]^n_+ }{ n! } d \beta(u),
\end{equation}
where $\beta(u)$ is non-decreasing.

We point out that the representation (\ref{eq:3.4}) has a short proof. Also, without loss of generality we may assume that $a=-\infty$ and $b=\infty$.
\begin{theorem}\label{th:3.0}
Let $n \geqslant 1$. A function $f \colon \real \to \real$ satisfying the condition $f(-\infty)=0$ is $n$-times monotone non-decreasing if and only if it admits the representation
\begin{equation}
f(x) = \int_{-\infty}^\infty \frac{(x-u)^{n-1}_+}{(n-1)!} d \beta(u),
\label{eq:3.5.a}
\end{equation}
where $\beta(u)$ is non-decreasing and $\beta(-\infty) = 0$. Moreover, $\beta(u)$ is unique at its points of continuity and $\beta(u) = f^{(n-1)}(u)$ a.e.
\end{theorem}

\begin{proof}
The sufficiency is evident by differentiating (\ref{eq:3.5.a}), for
$$
f^{(k)}(x) = \int_{-\infty}^\infty \frac{(x-u)^{n-1-k}_+}{(n-1-k)!} d\beta(u),
$$
$(k=0,1,\ldots,n-1)$ is evidently non-negative and non-decreasing.

To see the necessity let us consider $\beta(u) = f^{(n-1)}(u)$ (with the convention $f^{(0)}(x) = f(x)$). Then (\ref{eq:3.5.a}) can be rewritten as
\begin{equation}
f(x) = \int_{-\infty}^x \frac{(x-u)^{n-1}}{(n-1)!} df^{(n-1)}(u).
\label{eq:3.5b}
\end{equation}
We prove it by induction.

Let $n=1$. Since $f(x)$ is non-decreasing and $f(-\infty)=0$, we have
\begin{equation}
f(x) = \int_{-\infty}^x d f(u).
\label{eq:3.5c}
\end{equation}
This proves (\ref{eq:3.5b}) for $n=1$.

Assume that (\ref{eq:3.5b}) holds for some $n \geqslant 1$, i.e. $n$-times monotone non-decreasing functions $f(x)$ such that $f(-\infty)=0$ are of the form (\ref{eq:3.5b}). Let $f(x)$ (such that $f(-\infty)=0$) be an $(n+1)$-times monotone non-decreasing, i.e. $f(x)$ is $n$-times monotone non-decreasing and $f^{(n)}(x)$ is non-decreasing. It is not difficult to show that $f^{(n)}(-\infty)=0$. By (\ref{eq:3.5b}) and (\ref{eq:3.5c}), with $f^{(n)}(x)$ in place of $f(x)$, we have
\begin{eqnarray*}
f(x) &
= & \int_{-\infty}^x \frac{(x-u)^{n-1}}{(n-1)!} d f^{(n-1)}(u) \\
& = & \int_{-\infty}^x \frac{(x-u)^{n-1}}{(n-1)!} f^{(n)}(u) du \\
& = & \int_{-\infty}^x \frac{(x-u)^{n-1}}{(n-1)!} \int_{-\infty}^u d f^{(n)}(v) du \\
& = & \int_{-\infty}^x \int_{-\infty}^x \frac{(x-u)^{n-1}}{(n-1)!} du df^{(n)}(v) \\
& = & \int_{-\infty}^x \frac{(x-v)^n}{n!} df^{(n)}(v).
\end{eqnarray*}
This proves (\ref{eq:3.5b}), with $n$ replaced by $n+1$, so the induction is complete.

\end{proof}

Let $\Mcl_{n+}((a,b))$ ($\Mcl_{n-}((a,b))$) be the class of all $n$-times monotone non-decreasing (non-increasing) functions on $(a,b)$. Taking into account (\ref{eq:3.4}) and (\ref{eq:3.5}), by Theorem (\ref{th:2.9}) we obtain the following decomposition.

\begin{theorem}\label{th:3.1}
Let $n \geqslant 1$ and $f \colon (a,b) \to \real$. Then $f$ is $n$th order convex if and only if $f$ is of the form
$$ %\begin{equation}\label{eq:3.6}
f(x) = M_1(x) + M_2(x) + Q(x),
$$ %\end{equation}
where $(-1)^{n+1} M_1(x) \in \Mcl_{(n+1)-}((a,\xi))$, $M_2(x) \in \Mcl_{(n+1)+}((\xi,b))$, with $a \leqslant \xi \leqslant b$,
$Q(x) = c_n x^n/n!+\ldots+c_0 \in \Pi_n$. Moreover, 
if $M_1 \equiv 0$ and $M_2 \not\equiv 0$ then $c_n \geqslant 0 $,
if $M_2 \equiv 0$ and $M_1 \not\equiv 0$ then $c_n \leqslant 0 $, and
if $M_1 \not\equiv 0$ and $M_2 \not\equiv 0$ then $c_n = 0 $ and $M_1(\xi-) = M_2(\xi+)$.
\end{theorem}

%\begin{remark}\label{remark:3.1.a}
%Note that $M_1(x)$ is $n$-convex non-increasing for $n$ odd and non-decreasing for $n$ even, $M_2(x)$ is $n$-convex non-decreasing.
%\end{remark}

\begin{remark}\label{remark:3.1.1}
Note that if $g(x) \in \Mcl_{(n+1)+}((a,b))$, then $\varphi(x) = g(-x) \in \Mcl_{(n+1)-}((-b,-a))$. Thus $M_2(x) = g(x)$ is $n$th order convex on $(a,b)$ and $M_1(x) = (-1)^{n+1} \varphi(x)$ is $n$th order convex on $(-b,-a)$.
\end{remark}

The following proposition gives a characterization of $n$-times monotone non-decreasing functions in terms of difference operators (see McNeil \cite{McNeilNeslehova2009}).

\begin{proposition}\label{prop:3.2}
Let $n \geqslant 1$ and $f \colon (a,b) \to \real$. Then the following statements are equivalent.
\begin{itemize}
	\item[(i)] $f$ is $n$-times monotone on $(a,b)$.
	\item[(ii)] $f$ is non-negative and for any $k=1,\ldots,n$, any $x \in (a,b)$, and any $h_i > 0$, $i=1,\ldots,k$ such that $x+h_1+\ldots+h_k \in (a,b)$ the function $f$ satisfies
	\begin{equation}\label{eq:3.7}
	\Delta_{h_k}\ldots \Delta_{h_1} f(x) \geqslant 0,
	\end{equation}
	where $\Delta_{h_k}\ldots \Delta_{h_1}$ denote sequential applications of the first-order difference operator $\Delta_h$ given by $\Delta_h f(x) = f(x+h)-f(x)$ whenever $x, x+h \in (a,b)$.
\item[(iii)] $f$ is non-negative and satisfies, for any $k=1,\ldots,n$, any $x \in (a,b)$, and any $h>0$ such that $x+kh \in (a,b)$
 \begin{equation}\label{eq:3.8}
 (\Delta_h)^k f(x) \geqslant 0,
 \end{equation}
 where $(\Delta_h)^k$ denote the $k$-monotone sequential iterations of the operator $\Delta_h$.
\end{itemize}
\end{proposition}

Note that in Gilanyi and Pales \cite{GilanyiPales2008} functions satisfying (\ref{eq:3.7}) with $k=n+1$ are called \textit{Wright-convex of order $n$} (or simply \textit{$n$-Wright-convex}). As it is extensively discussed in \cite{Kuczma1985} ($n \geqslant 1$), the functions satisfying (\ref{eq:3.8}) with $k=n+1$ are called \textit{Jensen-convex of order $n$}. It is well-known that, under the assumption of continuity, Jensen convexity of order $n$ and $n$th order convexity are equivalent. In the study of inequalities (\ref{eq:3.7}) and (\ref{eq:3.8}), functions that satisfy (\ref{eq:3.7}) and (\ref{eq:3.8}) with equality play a crucial role. For $n \in \nat$, a function $P \colon \real \to \real$ is called a \textit{polynomial function of degree at most $n$} if it satisfies the Fr\'{e}chet equation, i.e. if
$$
(\Delta_h)^{n+1} P(x) = 0 \quad (h,x\in\real).
$$
Polynomials are exactly the continuous polynomial functions, however, in terms of Hamel bases, one can construct non-continuous polynomial functions (see \cite{Kuczma1985}). Maksa and P\'{a}les \cite{MaksaPales2009} proved that any $n$-Wright-convex function can be represented as the sum of a continuous $n$-convex function and a polynomial function
\begin{proposition}\label{prop:3.3}
Let $n \geqslant 1$ and $f \colon (a,b) \to \real$. Then $f$ is an $n$-Wright-convex function if and only if $f$ is of the form
\begin{equation}\label{eq:3.9}
f(x) = C(x) + P(x) \quad (x \in (a,b)),
\end{equation}
where $C \colon (a,b) \to \real$ is a continuous $n$-convex function and $P \colon \real \to \real$ is a polynomial function of degree at most $n$ with $P(\Qb) = \{0 \}$. Furthermore, under the assumption $P(\Qb) = 0$, the decomposition (\ref{eq:3.9}) is unique.
\end{proposition}

From Theorem \ref{th:3.1} and Proposition \ref{prop:3.3} we obtain the following decomposition of $n$-Wright-convex functions.

\begin{theorem}\label{th:3.4}
Let $n \geqslant 1$ and $f \colon (a,b) \to \real$. Then $f$ is an $n$-Wright-convex if and only if $f$ is of the form
$$
f(x) = M_1(x) + M_2(x) + Q(x) + P(x),
$$
where $(-1)^{n+1} M_1(x)$ is an $(n+1)$-times monotone non-increasing on $(a,\xi)$, $M_2(x)$ is an $(n+1)$-times monotone non-decreasing on $(\xi,b)$, $a \leqslant \xi \leqslant b$, $M_1(x) + M_2(x)$ is continuous on $(a,b)$, $Q(x)$ is a polynomial of degree at most $n$ (as in Theorem \ref{th:3.1}) and $P(x)$ is a polynomial function of degree at most $n$.
\end{theorem}

\section{Relative $n$-convexity. Strong $n$-convexity.}

Let $g \colon (a,b) \to \real$ be an $n$-convex function. We say that a function $f \colon (a,b) \to \real$ is \textit{$n$-convex with respect to $g$} if $f-g$ is $n$-convex, and denote it by $f \succeq_n g$.

Various other generalizations of convexity via related convexity properties have been proposed. The relative $n$-convexity defined above is a generalization of the relative convexity (for $n=1$) studied in Karlin and Studden \cite{KarlinStudden1966} (cf. \cite{Cargo1965}, \cite{HardyLittlewoodPolya1959}, \cite{Palmer2003}, \cite{PecaricProschanTong1992}).

\begin{remark}\label{remark:4a.1}
If $f$ is $n$-convex with respect to $g$, then both $f-g$ and $g$ are $n$-convex. Writting $f=g+(f-g)$, we obtain that $f$ necessarily must be $n$-convex.
\end{remark}

Functions $f \colon (a,b) \to \real$ and $g \colon (a,b) \to \real$ that are decreasing and increasing on the same intervals will be called \textit{isotonic}, and we say that they are members of the same isotonic class.

\begin{theorem}\label{th:4a.3}
Let $f,g \colon (a,b) \to \real$ be $n$-convex functions with $n$-spectral measures $\mu_{(n)-}$, $\mu_{(n)+}$ and $\nu_{(n)-}$, $\nu_{(n)+}$, respectively. Then $f$ is $n$-convex with respect to $g$ if and only if
\begin{equation}
\mu_{(n)} \geqslant \nu_{(n)},
\label{eq:4a.4}
\end{equation}
where
$$ %\begin{equation}
\mu_{(n)} = \mu_{(n)-} + \mu_{(n)+},
$$ %\label{eq:4a.5}
%\end{equation}
$$ %\begin{equation}
\nu_{(n)} = \nu_{(n)-} + \nu_{(n)+}.
%\label{eq:4a.6}
$$ %\end{equation}
\end{theorem}

\begin{proof}
Let $f$ and $g$ satisfy the assumptions of the theorem. Fix $a < \xi <  b$. By Theorem \ref{th:2.10}, $f(x)$ and $g(x)$ can be written in the form (\ref{eq:2.15}) with the measures $\mu_{(n)\xi-}$, $\mu_{(n)\xi+}$, $\nu_{(n)\xi-}$ and $\nu_{(n)\xi+}$, and the polynomials $Q_\xi(x)$ and $R_\xi(x)$, respectively. In other words, the functions $f(x) - Q_\xi(x)$ and $g(x) - R_\xi(x)$ are isotonic. Moreover, by (\ref{eq:2.17}), we have
\begin{equation}
\mu_{(n)-} + \mu_{(n)+} = \mu_{(n)\xi-} + \mu_{(n)\xi+}, \quad
\nu_{(n)-} + \nu_{(n)+} = \nu_{(n)\xi-} + \nu_{(n)\xi+}.
\label{eq:4a.7}
\end{equation}
Therefore $f(x)-g(x)$ is of the form (\ref{eq:2.15}), with the measures $\mu_{(n)\xi-} - \nu_{(n)\xi-}$ and $\mu_{(n)\xi+} - \nu_{(n)\xi+}$ in the place of $\mu_{(n)\xi-}$ and $\mu_{(n)\xi+}$, respectively, and $Q_\xi(x) - R_\xi(x)$ in the place of $Q_\xi(x)$. By Theorem \ref{th:2.10}, $f-g$ is $n$-convex if and only if
$$
\mu_{(n)\xi-} \geqslant \nu_{(n)\xi-}, \quad 
\mu_{(n)\xi+} \geqslant \nu_{(n)\xi+}.
$$
From (\ref{eq:4a.7}) we conclude (\ref{eq:4a.4}). The theorem is proved.
\end{proof}

Theorem \ref{th:4a.3} suggests that we can define a \textit{measure of $n$th order convexity} of an $n$-convex function by the operator
$$
K \colon f \to \mu_{(n)}^f = \mu_{(n)}=\mu_{(n)-}+\mu_{(n)+}.
$$
In the sequel we will call $\mu_{(n)}^f$ the \textit{measure of $n$-convexity of $f$}, or shortly the \textit{$n$-convexity measure}. From Theorem \ref{th:4a.3} we have
\begin{theorem}\label{th:4.2.1}
$$
f \succeq_n g \textit{ if and only if } \mu_{(n)}^f \geqslant \mu_{(n)}^g.
$$
\end{theorem}

We shall say that functions $f,g \colon (a,b) \to \real$ are of \textit{modulo $\Pi_n$}, or that they are members of the same \textit{modulo $\Pi_n$ class}, if they differ by a polynomial $Q \in \Pi_n$. The relation modulo $\Pi_n$ is an equivalence relation and hence it defines equivalence classes. For $n$-convex $f$ and $g \colon (a,b) \to \real$ that are members of the same modulo $\Pi_n$ class we therefore have that $f^{(n)}(x)$ and $g^{(n)}(x)$ differ on $(a,b)$ by a constant. Consequently, by Theorem \ref{th:2.10}, we have the following theorem

\begin{theorem}\label{th:4a.4}
$$
f=g \; (\textrm{mod}\;\Pi_n) \quad \textit{ if and only if } \quad \mu_{(n)}^f = \mu_{(n)}^g.
$$
\end{theorem}
We now show that this relation induces a partial ordering.
\begin{theorem}\label{th:4a.5}
The relative $n$-convexity relation induces a partial ordering on modulo $\Pi_n$ equivalence classes of $n$-convex functions.
\end{theorem}

\begin{proof}
We will show that the relation is reflective, antisymmetric, and transitive.

\noindent
\textbf{Reflectivity.} For all $f$ we have $f-f \equiv 0 \in \Pi_n$. Thus $f \succeq_n f$.

\noindent
\textbf{Antisymmetry.} Suppose $f \succeq_n g$ and $g \succeq_n f$. Then $f-g$ and $g-f$ are $n$-convex. Thus, both functions $(f-g)^{(n)}(x)$ and $[-(f-g)]^{(n)}(x)$ are non-decreasing. Consequently, $(f-g)^{(n)}(x) = 0$ ($x \in (a,b)$). This implies that $f-g \in \Pi_n$, that is $f=g \; (\textrm{mod}\;\Pi_n)$.

\noindent
\textbf{Transitivity.} Suppose $f \succeq_n g$ and $g \succeq_n h$. Then both $f-g$ and $g-h$ are $n$-convex. Writing $f-h$ in the form $f-h = (f-g)+(g-h)$, we obtain that $f-h$ is $n$-convex as the sum of the $n$-convex functions. The theorem is proved.
\end{proof}

As a simple example of the use of Theorem \ref{th:2.10} we prove the following theorem. We denote by $d \mu/d \nu = \varphi$ the Radon-Nikodym derivative of a measure $\mu$ with respect to a measure $\nu$ (see \cite{Royden1966}).

\begin{theorem}\label{th:5.5a}
Let $f,g \colon (a,b)\to\real$ be $n$-convex. Then
\begin{enumerate}[\upshape a)]
	\item there exists an $n$-convex function $f_{\textit{max}}$ such that
	$$
	f_{\textit{max}} \succeq_n f, \quad f_{\textit{max}} \succeq_n g,
	$$
	and for every $n$-convex function $h$
	$$
	( h \succeq_n f \textit{ and } h \succeq_n g ) \Rightarrow h \succeq_n f_{\textit{max}},
	$$

	\item there exists an $n$-convex function $f_{\textit{min}}$ such that
	$$
	f \succeq_n f_{\textit{min}}, \quad g \succeq_n f_{\textit{min}},
	$$
	and for every $n$-convex function $h$
	$$
	( f \succeq_n h \textit{ and } g \succeq_n h ) \Rightarrow  f_{\textit{min}} \succeq_n h,
	$$

	\item if $f \succeq_n g$ and $f \neq g (\textrm{mod}\;\Pi_n)$, then there exists an $n$-convex function $w$ such that $f \neq w (\textrm{mod}\;\Pi_n)$, $g \neq w (\textrm{mod}\;\Pi_n)$ and 
	$$
	f \succeq_n w \succeq_n g.
	$$
\end{enumerate}
\end{theorem}

\begin{proof}
Let $f$ and $g$ be $n$-convex. By Theorem \ref{th:2.10} we may assume that $f$ and $g$ admit representations given by (\ref{eq:2.15}) with the same $\xi \in (a,b)$, the measures $\mu_{(n)\xi-}^f$, $\mu_{(n)\xi+}^f$, $\mu_{(n)\xi-}^g$, $\mu_{(n)\xi+}^g$, and with the polynomials $Q_\xi^f$ and $Q_\xi^g$, respectively. Consider the Radon-Nikodym derivatives
$$
\varphi_1 = d \mu_{(n)\xi-}^f / d( \mu_{(n)\xi-}^f + \mu_{(n)\xi-}^g ),
$$
$$
\psi_1 = d \mu_{(n)\xi-}^g / d( \mu_{(n)\xi-}^f + \mu_{(n)\xi-}^g ),
$$
$$
\varphi_2 = d \mu_{(n)\xi+}^f / d( \mu_{(n)\xi+}^f + \mu_{(n)\xi+}^g ),
$$
$$
\psi_2 = d \mu_{(n)\xi+}^g / d( \mu_{(n)\xi+}^f + \mu_{(n)\xi+}^g ).
$$
It is not difficult to see that it suffices to take the functions $f_{\textit{max}}$ and $f_{\textit{min}}$ of the form (\ref{eq:2.15}) with the measures
$$
\mu_{(n)\xi-}^{\textit{max}} = \max( \varphi_1,\psi_1 ) ( \mu_{(n)\xi-}^{f} + \mu_{(n)\xi-}^{g} ),
$$
$$
\mu_{(n)\xi+}^{\textit{max}} = \max( \varphi_2,\psi_2 ) ( \mu_{(n)\xi+}^{f} + \mu_{(n)\xi+}^{g} ),
$$
$$
\mu_{(n)\xi-}^{\textit{min}} = \max( \varphi_1,\psi_1 ) ( \mu_{(n)\xi-}^{f} + \mu_{(n)\xi-}^{g} ),
$$
$$
\mu_{(n)\xi+}^{\textit{min}} = \max( \varphi_2,\psi_2 ) ( \mu_{(n)\xi+}^{f} + \mu_{(n)\xi+}^{g} ),
$$
to prove parts a) and b).

To prove c) assume $f \succeq_n g$ and $f \neq g (\textrm{mod}\; \Pi_n)$. Then $f-g$ is $n$-convex and $f-g \neq 0 (\textrm{mod}\;\Pi_n)$. Thus it suffices to take $w = g+\frac{1}{2}(f-g)$. The theorem is proved.
\end{proof}

As usual we denote distributional derivatives by $f'$ (see \cite{Schwartz1966}, \cite{Talvila2009}), pointwise derivatives by $f'(x)$, $n$th order distributional derivatives by $f^{(n)}$, and $n$th order pointwise derivatives by $f^{(n)}(x)$. Theorem (\ref{th:4a.3}) suggests that we can use distributional derivatives and the Radon-Nikodym derivatives to derive simple criteria for the relative $n$-convexity $f \succeq_n g$.

\begin{theorem}\label{th:4a.6}
Let $f,g \colon (a,b) \to \real$ be $n$-convex functions with the $n$-convexity measures $\mu_{(n)}^f$ and $\mu_{(n)}^g$, respectively. Then the following conditions are equivalent:
\begin{enumerate}[\upshape a)]
	\item $f \succeq_n g$,
	\item $\mu_{(n)}^f \geqslant \mu_{(n)}^g$,
	\item $f^{(n+1)} \geqslant g^{(n+1)}$,
	\item $d\mu_{(n)}^g / d\mu_{(n)}^f \leqslant 1$.
\end{enumerate}
\end{theorem}

Via Lebesgue's decomposition theorem and the decomposition of a singular measure, every $\sigma$-finite measure $\mu$ can be decomposed into a sum of an absolutely continuous measure (with respect to the Lebesgue measure), a singular continuous measure, and a discrete measure, i.e.
$$
\mu = \mu_{\textit{cont}} + \mu_{\textit{sing}} + \mu_{\textit{pp}},
$$
where $\mu_{\textit{cont}}$ is the absolutely continuous part, $\mu_{\textit{sing}}$ is the singular continuous part and $\mu_{\textit{pp}}$ is the pure point part (a discrete measure) (see Royden \cite{Royden1966}). These three measures are uniquely determined.

\begin{remark}\label{remark:5.6a}
The following decomposition yields an analoguous decomposition of an $n$-convex function. Namely, any $n$-convex function $f$ with the $n$-spectral measures $\mu_{(n)-}$ and $\mu_{(n)+}$ can be represented as a sum 
\begin{equation}
f = f_{\textit{cont}} + f_{\textit{sing}} + f_{\textit{pp}},
\label{eq:4a.8}
\end{equation}
where $f_{\textit{cont}}$, $f_{\textit{sing}}$ and $f_{\textit{pp}}$ correspond to the absolutely continuous parts the singular continuous parts, and the pure point parts of the $n$-spectral measures $\mu_{(n)-}$ and $\mu_{(n)+}$, respectively 
($\mu_{(n)-}=\mu_{(n)-\textit{cont}} + \mu_{(n)-\textit{sing}} + \mu_{(n)-\textit{pp}}$,
$\mu_{(n)+}=\mu_{(n)+\textit{cont}} + \mu_{(n)+\textit{sing}} + \mu_{(n)+\textit{pp}}$).
Note, that $f_{\textit{cont}}$, $f_{\textit{sing}}$ and $f_{\textit{pp}}$ are $n$-convex. Moreover, they are unique, up to a polynomial of degree at most $n$.
\end{remark}

It is not difficult to prove the following lemma.

\begin{lemma}\label{lemma:4a.6}
Let $\mu$ and $\nu$ be two $\sigma$-finite measures having the following decompositions into a sum of an absolutely continuous measure, a singular continuous measure and a discrete measure
$$
\mu = \mu_{\textit{cont}} + \mu_{\textit{sing}} + \mu_{\textit{pp}} \quad \textit{and }\quad
\nu = \nu_{\textit{cont}} + \nu_{\textit{sing}} + \nu_{\textit{pp}}.
$$
Then $\mu \geqslant \nu$ if and only if $\mu_{\textit{cont}} \geqslant \nu_{\textit{cont}}$, $\mu_{\textit{sing}} \geqslant \nu_{\textit{sing}}$ and $\mu_{\textit{pp}} \geqslant \nu_{\textit{pp}}$.
\end{lemma}

Taking into account the decomposition (\ref{eq:4a.8}) of an $n$-convex function, by Lemma \ref{lemma:4a.6} we immediately obtain the following three theorems useful in studying relative $n$-convexity.

\begin{theorem}\label{th:4a.7}
Let $f$ and $g \colon (a,b) \to \real$ be $n$-convex functions having the decompositions $f=\fcont+\fsing+\fpp$, $g=\gcont+\gsing+\gpp$ (see (\ref{eq:4a.8})). Then $f \succeq_n g$ if and only if $\fcont \succeq_n \gcont$ and $\fsing \succeq_n \gsing$ and $\fpp \succeq_n \gpp$.
\end{theorem}

\begin{theorem}\label{th:4a.8}
$$
\fcont \succeq_n \gcont \quad \textit{ iff } \quad
\fcont^{(n+1)}(x) \geqslant \gcont^{(n+1)}(x)
$$
\end{theorem}

\begin{theorem}\label{th:4a.9}
$$
\fpp \sn \gpp \quad \textit{ iff } \quad \fpp^{(n+1)} \geqslant \gpp^{(n+1)},$$
where $\fpp^{(n+1)} = \sum_{k} a_k \delta_{x_k}$, 
$\gpp^{(n+1)} = \sum_{k} b_k \delta_{y_k}.$
\end{theorem}

The notion of convexity can be extended not only to the case when the order of convexity is of higher-dimension, but also in several other ways. One of the most important generalizations is the notion of strong convexity. A function $f \colon (a,b) \to \real$ is called \textit{strongly convex with modulus $c > 0$} if
$$
f( tx + (1-t)y ) \leqslant t f(x) + (1-t) f(y) - ct(1-t)(x-y)^2
$$
for all $x,y \in (a,b)$ and $t \in [0,1]$. Strongly convex functions	were introduced by Polyak \cite{Polyak1966}. Some properties of them can be found, among other, in \cite{RobertsVarberg1973}, \cite{HiriartLemarechal2001}, \cite{Polovinkin1966}. Not attempting to be complete, we just recall here two results concerning strong convexity which play crucial roles in further invesigations (see \cite{RobertsVarberg1973}). The first one characterizes strong convexity in terms of convexity, while the second one characterizes twice differentiable strongly convex function in terms of its second derivative $f''(x)$.

\begin{proposition}\label{prop:5.10a}
A function $f \colon (a,b) \to \real$ is strongly convex with modulus $c > 0$ if and only if the function $f(x) - cx^2$ is convex.
\end{proposition}

\begin{proposition}\label{prop:5.10b}
Assume that $f \colon (a,b) \to \real$ is twice differentiable and $c>0$. Then $f$ is strongly convex with modulus $c$ if and only if $f''(x) \geqslant 2c$ ($x \in (a,b)$).
\end{proposition}

As a generalization of strong convexity with modulus $c$, we define strong $n$-convexity with modulus $c$. We say that a function $f$ is \textit{strongly $n$-convex with modulus $c$} ($n \geqslant 1$, $c>0$) if $f$ is $n$-convex with respect to the function $g(x) = c x^{(n+1)}/(n+1)!$. By Proposition \ref{prop:5.10a} the strong convexity with modulus $2c$ (cf Roberts and Varberg \cite{RobertsVarberg1973}) coincides with our strong $1$-convexity with modulus $c$. Writing $f(x) = \left( f(x) - c x^{n+1}/(n+1)! \right) + cx^{n+1}/(n+1)!$, we obtain that if $f$ is strongly $n$-convex with modulus $c>0$, then $f$ is $n$-convex.

The following theorem gives a characterization of a strongly $n$-convex function $f$ with modulus $c$ without additional assumptions on differentiability of $f$. This generalizes Proposition \ref{prop:5.10b}.

\begin{theorem}\label{th:4a.11}%Th. 5.12
Let $f \colon (a,b) \to \real$ be an $n$-convex function and $c>0$. Then $f$ is strongly $n$-convex with modulus $c$ if and only if
$$
f^{(n+1)}(x) \geqslant c \textit{ for } x \in (a,b) \quad \lambda\; a.e.
$$
\end{theorem}

\begin{proof}
Note that the function $g(x) = c x^{n+1}/(n+1)!$ is $n$-convex with the $n$-convexity measure $\mu^g_{(n)} (dx) = dg^{(n)}(x) = cdx$, where $g^{(n)}(x) = cx$. Writting $g(x)$ in the form (\ref{eq:4a.8}), we have $g = g_{\textit{cont}}$, $g_{\textit{sing}} = 0$ and $g_{\textit{pp}} = 0$. As $f$ is $n$-convex, we look at its integral representation given by (
\ref{eq:2.9}), with the $n$-convexity measure $\mu_{(n)}^f (dx) = df^{(n)}(x)$. Since $f^{(n)}(x)$ is nondecreasing its derivative $f^{(n+1)}(x)$ exists for $x \in (a,b)$ $\lambda$ a.e. By Remark \ref{remark:5.6a} $f(x)$ can be represented as the sum
$$
f = f_{\textit{cont}} + f_{\textit{sing}} + f_{\textit{pp}}.
$$
From Theorem \ref{lemma:4a.6}, and taking into account that $g_{\textit{cont}} = g$, $g_{\textit{sing}} = g_{\textit{pp}} = 0$, we obtain that
$$
f \succeq_n g \textit{ iff } f_{\textit{cont}} \succeq_n g_{\textit{cont}}.
$$
By Theorem \ref{th:4a.7}, $f_{\textit{cont}} \succeq_n g_{\textit{cont}}$ iff $f_{\textit{cont}}^{(n+1)}(x) \geqslant g_{\textit{cont}}^{(n+1)}(x)$. Since $g_{\textit{cont}}^{(n+1)}(x) = c$, and $f_{\textit{cont}}^{(n+1)}(x) = f^{(n+1)}(x)$ for $x \in (a,b)$ $\lambda$ a.e., the theorem is proved.
\end{proof}

\begin{corollary}\label{corollary:5.13}
Let $c > 0$, $n \in \nat$ and $f \colon (a,b) \to \real$ be a function. Then $f$ is strongly $n$-convex with modulus $c$ if and only if $f$ is of the form
$$
f(x) = f_{\textit{cont}}(x) + R(x) \quad (x \in (a,b)),
$$
where $f_{\textit{cont}} \colon (a,b) \to \real$ is an $(n+1)$-times differentiable strongly $n$-convex function with modulus $c$, and $R \colon (a,b) \to \real$ is an $n$-convex function such that $R^{(n+1)}(x) = 0$ for $x \in (a,b)$ $\lambda$ a.e.
\end{corollary}

\begin{corollary}\label{corollary:5.14}
Let $c>0$, $n\in\nat$ and let $f \colon (a,b) \to \real$ be an $(n+1)$-times differentiable function. Then $f$ is strongly $n$-convex with modulus $c$ if and only if
$$
f^{(n+1)}(x) \geqslant c, \quad x \in (a,b).
$$
\end{corollary}

\section{Interpolation of functions by $n$-convex functions.}

It is well-known that every convex function $f \colon I \to \real$ admits an affine support at every interior point of $I$ (i.e. for any $x_0 \in \textit{Int} I$ there exists an affine function $a \colon I \to \real$ such that $a(x_0) = f(x_0)$ and $a \leqslant f$ on $I$). Convex functions of higher orders (precisely of an odd orders) have similar property; they are supported by polynomials of degree no greater than the order of convexity.

The following important property of convex functions of higher order (cf. Kuczma \cite{Kuczma1985}, Popoviciu \cite{Popoviciu1934}, Roberts and Varberg \cite{RobertsVarberg1973}) is well-known: a function $f \colon I \to \real$ is $n$-convex ($I \subset \real$ is an interval) if and only if for any $x_1,\ldots,x_{n+1} \in I$ with $x_1 < \ldots x_{n+1}$ the graph of an interpolating polynomial $p := P(x_1,\ldots,x_{n+1};f)$ passing through the points $(x_i, f(x_i))$, $i=1,\ldots,n+1$, changes successively from one side of the graph of $f$ to another (always $p(x) \leqslant f(x)$ for $x \in I$ such that $x > x_{n+1}$ if such points exist). More precisely,
$(-1)^{n+1} (f(x) - p(x) \geqslant 0)$ $x < x_1$, $x \in I$, $(-1)^{n+1-i} ( f(x) - p(x) ) \geqslant 0$, $x_i < x < x_{i+1}$, $i=1,\ldots,n$,
$f(x)- p(x) \geqslant 0$, $x > x_{n+1}$, $x\in I$. It is not difficult to observe that the $n$-convexity reduces to convexity in the usual sense if $n=1$.

In the Wasowicz paper \cite{Wasowicz2007} certain attaching method is developed. The method is applied in Theorem \ref{th:5.3} to obtain a general result, from which the mentioned above support theorem and some related properties of convex functions of higher order are derived.
\begin{theorem}\label{th:5.3}
Let $n \in \nat$ and $f \colon I \to \real$ be an $n$-convex function. Let us fix $k \in \nat$, $k \leqslant n$, and take $x_1,\ldots,x_k \in I$ such that $x_1 < \ldots < x_k$. To each point $x_j$ $(j=1,\ldots,k)$ assign the multiplicity $l_j \in \nat$ such that $l_1 + \ldots l_j = n+1$. Additionally assume that if $x_1 = \inf I$, then $l_1=1$, and if $x_k = \sup I$, then $l_k=1$. Denote $I_0 = (-\infty,x_1)$, $I_j = (x_j,x_{j+1})$, $j=1,\ldots,k-1$, and $I_k=(x_k,\infty)$. Under these assumptions there exists a polynomial $p \in \Pi_n$ such that $p(x_j) = f(x_j)$, $j=1,\ldots,k$, and such that 
$(-1)^{n+1} ( f(x) - p(x) ) \geqslant 0$ for $x \in I_0 \cap I$,
$(-1)^{n+1-(l_1 + \ldots l_j)}( f(x) - p(x) ) \geqslant 0$ for $x \in I_j$, $j=1,\ldots,k-1$,
$f(x) - p(x) \geqslant 0$ for $x \in I_k \cap I$.
\end{theorem}

The numbers $l_1,\ldots,l_k$ can be interpreted as multiplicities of the points $x_1,\ldots,x_k$, respectively. The polynomial $p(x)$ in the above theorem will be called the \textit{support of $(l_1,\ldots,l_k)$-type}.

\begin{remark}\label{remark:5.4}
This fact is shown by Wasowicz \cite{Wasowicz2010} in a more general setting, i.e. for functions convex with respect to Chebyshev systems (for Chebyshev's polynomial system $(1,x,\ldots,x^n)$ such convexity reduces to $n$-convexity).
\end{remark}

\begin{observation}\label{obs:5.5}
The polynomial $p(x)$ described in Theorem \ref{th:5.3} has following properties:
\begin{itemize}
	\item[(i)] $p(x) \leqslant f(x)$, $x > x_k$, $x \in I$,
	\item[(ii)] if $l_j$ (i.e. the multiplicity of $x_j$) is even, then the graph of $p(x)$ passing through $x_j$ remains on the same side of the graph of $f$, while it changes the side, if $l_j$ is odd.
\end{itemize}
\end{observation}

We apply Theorem \ref{th:5.3} to obtain a general result, that for any two $n$-convex functions $f$ and $g$, such that $f$ is $n$-convex with respect to $g$, the function $g$ is a support of $(l_1,\ldots,l_k)$-type for function $f$, up to some polynomial $p \in \Pi_n$.

\begin{theorem}\label{th:5.6}
Let $n \in \nat$ and let $f$ and $g \colon I \to \real$ be two $n$-convex functions such that $f$ is $n$-convex with respect to $g$. Fix $k \in \nat$, $k \leqslant n$, and let $x_1,\ldots,x_k \in I$ be such that $x_1 < \ldots < x_k$. Suppose that $l_j$, $I_j$ satisfy conditions of Theorem \ref{th:5.3}. Then there exists a polynomial $p \in \Pi_n$, such that
\begin{equation}
	f( x_j ) = g( x_j ) + p (x_j ), \quad j=1,\ldots,k,
	\label{eq:5.7}
\end{equation}
and additionally
$$
(-1)^{n+1} [ f(x) - (g(x) + p(x) )] \geqslant 0 \quad \textit{ for } x \in I_0 \cap I
$$
\begin{equation}
	(-1)^{n+1-(l_1+\ldots+l_j)}[ f(x) - (g(x)+p(x))] \geqslant 0 \quad \textit{ for } x \in I_j, j=1,\ldots,k-1,
	\label{eq:5.8}
\end{equation}
$$
f(x) - ( g(x) + p(x) ) \geqslant 0 \quad \textit{ for } x \in I_k \cap I.
$$
The function $g(x) + p(x)$ will be called the support of $(l_1,\ldots,l_k)$-type for the function $f$.
\end{theorem}

\begin{proof}
Since $f(x)$ is $n$-convex with respect to $g$, $f(x) - g(x)$ is $n$-convex. Applying Theorem \ref{th:5.3} with the function $f(x)-g(x)$ in place of $f(x)$, we obtain that there exsits a polynomial $p \in \Pi_n$ such that
$$
f(x_j)-g(x_j) = p(x_j) \quad j=1,\ldots,k,
$$
$$
(-1)^{n+1}[ (f(x)-g(x))-p(x)]\geqslant 0 \quad \textit{ for } x \in I_0 \cap I,
$$ 
$$
(-1)^{n+1-(l_1+\ldots+l_j)} [ (f(x)-g(x))-p(x) ]\geqslant 0 \quad \textit{ for } x \in I_j, j=1,\ldots,k-1,
$$
$$
(f(x)-g(x))-p(x) \geqslant 0 \quad \textit{ for } x \in I_k \cap I.
$$
Thus (\ref{eq:5.7}) and (\ref{eq:5.8}) are satisfied. This completes the proof.
\end{proof}



%% References


\bibliographystyle{elsarticle-num}



\end{document}